%% file: autosam10.tex
\documentclass[twocolumn]{autart}    % Enable this line and disable the 
\usepackage{natbib}
\usepackage{graphicx}          % Include this line if your 
\input{myIncludes}

\allowdisplaybreaks

\begin{document}
	
\raggedbottom

\begin{frontmatter}

\title{A geometric obstruction to almost global synchronization on Riemannian manifolds\thanksref{footnoteinfo}} % Title, preferably not more 

\thanks[footnoteinfo]{Corresponding author Johan Markdahl. Tel. +352 46 66 44 5085. 
Fax +352 46 66 44 35085.}

\author[LCSB]{Johan Markdahl}\ead{markdahl@kth.se}    % Add the 

\address[LCSB]{Luxembourg Centre for Systems Biomedicine, University of Luxembourg, Belval, Luxembourg.}            

% The template generates 'key words' but in the published papers it says 'keywords'. Oh well.          
                             % keyword list or with the 
                                          % help of the Automatica 
                                          % keyword wizard

\begin{abstract}
Multi-agent systems on nonlinear spaces sometimes fail to synchronize. This is usually attributed to the initial configuration of the agents being too spread out, the graph topology having certain undesired symmetries, or both. Besides nonlinearity, the role played by the geometry and topology of the nonlinear space is often overlooked. This paper concerns two gradient descent flows of quadratic disagreement functions on general Riemannian manifolds. One system is intrinsic while the other is extrinsic. We derive necessary conditions for the agents to synchronize from almost all initial conditions when the graph used to model the network is connected. If a Riemannian manifold contains a closed curve of locally minimum length, then there is a connected graph and a dense set of initial conditions from which the intrinsic system fails to synchronize. The extrinsic system fails to synchronize if the manifold is multiply connected. The extrinsic system appears in the Kuramoto model on $\smash{\mathcal{S}^1}$, rigid-body attitude synchronization on $\mathsf{SO}(3)$, the Lohe model of quantum synchronization on the $n$-sphere, and the Lohe model on $\U$. Except for the Lohe model on the $n$-sphere where $n\in\N\backslash\{1\}$, there are dense sets of initial conditions on which these systems fail to synchronize. The reason for this difference is that the $n$-sphere is simply connected for all $n\in\N\backslash\{1\}$ whereas the other manifolds are multiply connected. %These results show that in addition to the key distinction between synchronization on linear and nonlinear manifolds, there is an important subdivision between simply connected and non-simply connected nonlinear manifolds.

% max 1600 characters
\end{abstract}

\end{frontmatter}

\section{Introduction}

The study of emergent behaviour in complex systems that evolve on nonlinear spaces requires a geometric theory of synchronization. This paper uses the language of Riemannian geometry to formulate consensus protocols as gradient descent flows on manifolds. We explore how the global stability properties of the synchronized state is affected by imposing assumptions on the geometry and topology of the manifolds. While previous research has emphasized the distinction between linear and nonlinear spaces \citep{sarlette2009consensus, tron2013riemannian}, it has not resulted in further categorization. Still, much has been revealed with respect to specific manifolds, in particular for special instances of the Stiefel manifold, see \eg \citep{canale2007gluing,markdahl2016cdc,deville2018synchronization}. This paper unifies such results by connecting them to the existence or non-existence of a closed curve of locally minimum length in the manifold. In particular, a manifold being multiply connected implies that such a curve exists. The geometry and topology of the manifold is preventing the agents from synchronizing. The results of this paper explains why synchronization is almost global on all spheres except the circle \citep{markdahl2016cdc} and why the Kuramoto model on the circle is multistable \citep{canale2007gluing}. It also explains why almost global attitude synchronization on $\SOT$ cannot be achieved using the most natural consensus protocol \citep{sarlette2009consensus}. More generally, the same argument applies to synchronization on $\SO$ \citep{deville2018synchronization} and synchronization on the unitary group $\mathsf{U}(n)$, a manifold that appears in the Lohe model of quantum synchronization \citep{lohe2010quantum,ha2018relaxation}.

This paper builds on previous work on almost global consensus on nonlinear spaces initiated by  \citet{scardovi2007synchronization,sarlette2009consensus}. While it is known that almost global synchronization on nonlinear spaces is graph dependent in general \citep{sepulchre2011consensus}, this property has not previously been connected with the geometry and topology of the manifold. Moreover, this paper also furthers the work initiated by \citet{tron2012intrinsic,tron2013riemannian} on intrinsic consensus on general Riemannian manifolds. Previously, local stability of the synchronized state has been established. In this paper, we show that if a Riemannian manifold contains a closed curve of locally minimum length then the system is multistable; it contains at least one other Lyapunov stable equilibrium set aside from consensus. To the best of the author's knowledge, this is the first result that characterizes a stable equilibrium set aside from consensus in the case of intrinsic synchronization on general Riemannian manifolds.

There is a geometric approach to multi-agent consensus which is based on notions of convexity  \citep{moreau2005stability}. Typical results pertain to the case when all agents belong to a geodesically convex set \citep{afsari2011riemannian,tron2013riemannian,hartley2013rotation,chen2014finite}. For example, convergence to the consensus manifold of the $n$-sphere from any open  hemisphere under various assumptions on the graph topology \citep{zhu2013synchronization,lageman2016consensus,thunberg2018lifting,zhang2018equilibria}. A key property of such results is that the convex hull of the states of all agents shrinks over time. This approach is however restricted by the fact that the convex hull is not always defined on a global level. For example, the largest geodesically convex set on a sphere is an open hemisphere. However, the synchronized state is actually almost globally stable \citep{markdahl2018tac}. More powerful tools are hence needed to further our understanding of synchronization on nonlinear spaces.

This paper provides an intuitive notion for thinking about the evolution of multi-agent systems on a global level. For $N$ agents connected by a cycle graph on $N$ nodes, it is helpful to imagine the agents as beads on a string. Technically, the string would consists of $N$ geodesic curves, each connecting a pair of neighboring agents. If the manifold is simply connected, then a continuous shortening of the string to a point results in consensus. If not, then the agents must be threaded to one end of the string. This requires two neighboring agents to move away from each other, which goes against the basic design principle of consensus protocols. Indeed, the canonical consensus protocols obtained by generalizing the Kuramoto model to  manifolds sometimes fail to synchronize. %To guarantee almost global convergence on multi connected mannifolds requires {\it ad hoc} control design \citep{sarlette2009geometry,tron2012intrinsic}. 
This notion also leads us to see a connection between consensus seeking system on cycle graphs and curve shortening flows, a topic that has been widely studied in mathematics \citep{white2002evolution}. From this perspective, multi-agent consensus can be considered as a kind of polygon shortening flow  \citep{smith2007curve}. There are however some subtle differences between curve shortening flows and multi-agent systems on cycle graphs. Only in the case of large $N$ and under certain initial conditions does a multi-agent system over a cycle graph behave as a piecewise smooth curve.

%\noindent This paper shows that the difference between the circle and higher-dimensional spheres, \ie that convergence is almost global on $\smash{\mathcal{S}^n}$ for all $\smash{n\in\N\backslash\{1\}}$ but not on $\smash{\mathcal{S}^1}$, can be explained by geometric and topological properties of the manifolds. If a manifold contains a curve of minimum length, then there is a graph for which there exists a stable equilibrium manifold $\mathcal{Q}$ such that $\mathcal{Q}\cap\mathcal{C}=\varnothing$. This implies that the consensus manifold $\mathcal{C}$ is not globally asymptotically stable. A special case of this condition occurs when the manifold is closed but not simply connected. A manifold is simply connected if any closed curve on it can be continuously deformed into a point. In particular, neither the circle $\smash{\mathsf{S}^1}$ nor $\SO$ are simply connected, but all spheres of dimension $n\geq2$ are. While the condition of simple connectedness is necessary for almost global synchronization for all connected graphs, it is not sufficient in general as we show by counterexample. Simulations explore the extent to which simple connectedness may be sufficient on Stiefel manifolds.

\section{Preliminaries}

\label{sec:problem}

Let $(\M,g)$ be a Riemannian manifold. The set $\M$ is a real, smooth manifold and the metric tensor $g_x$ is an inner product on the tangent space $\ts[\M]{x}$ at $x$. The map $x\mapsto g_x(X,Y)$ is smooth for any two differentiable vector fields $X,Y$ on $\M$. The shortest curve $\gamma:[a,b]\rightarrow\M$ such that $\gamma(a)=x$, $\gamma(b)=y$ is a geodesic from $x$ to $y$ (up to parametrization). The length of a curve $\gamma:[a,b]\rightarrow\M$ is
\begin{align}
l(\gamma)=\int_{a}^b g_\gamma(\dot{\gamma},\dot{\gamma})^\frac12\diff t.\label{eq:l}
\end{align}
Let $\Gamma$ denote the set of piece-wise smooth curves on $\M$. The length of a geodesic curve is the geodesic distance
\begin{align*}
d_g(x,y)=\inf\{l(\gamma)\,|\,\gamma\in\Gamma,\,\gamma(a)=x,\gamma(b)=y\}.
\end{align*}
When using the term geodesics in this paper we are only concerned with this property, and disregard the requirement of unit speed  parameterization. We extend the notion of geodesic distance to sets, defining
\begin{align*}
d_g(\mathcal{X},\mathcal{Y})=\inf\{l(\gamma)\,|\,\gamma\in\Gamma,\,\gamma(a)\in \mathcal{X},\gamma(b) \in \mathcal{Y}\}.
\end{align*}
Of particular interest are manifolds that have been embedded in the ambient space $\R^{n\times m}$. For such any manifold $\M\subset\R^{n\times}$, we use bold font to denote any elements $\ma{X},\ma{Y}\in\M$ and define the chordal distance $d_c(\ma{X},\ma{Y})=\|\ma{X}-\ma{Y}\|$ based on the Frobenius norm $\|\cdot\|:\ma{X}\mapsto\smash{(\trace\mat{X}\ma{X})^{\frac12}}$.

There is a Hilbert manifold of maps from the circle $\mathcal{S}^1$ to $\M$ \citep{klingenberg1978lectures}. That framework allows for a theory of closed geodesics on Riemannian manifolds. We will not detail it here, but provide the following result, adapted to the context of this paper:
\begin{theorem}[\citet{klingenberg1978lectures}]\label{th:simply}
Assume that the Riemannian manifold $(\M,g)$ is closed and multiply connected. Then $(\mathcal{M},g)$ contains a closed curve that is a local minimizer of the curve length function $l$ given by \eqref{eq:l}.
\end{theorem}
\noindent The property of simple connectedness refers to a path connected manifold on which each closed curve can be continuously deformed to a point. A multiply connected manifold is path-connected but contains at least one closed curve which cannot be continuously deformed to a point. In this paper, we focus on manifolds $\M$ that are path connected since almost global consensus would be impossible otherwise. If closedness is omitted from the requirements of Theorem \ref{th:simply}, then a counterexample is given by the punctured plane $ \smash{\R^2}\backslash\{\ve{0}\}$. This manifold is multiply connected, yet it does not contain a closed curve of locally minimum length.

\begin{figure}[htb!]
	\centering	\includegraphics[width=0.45\textwidth]{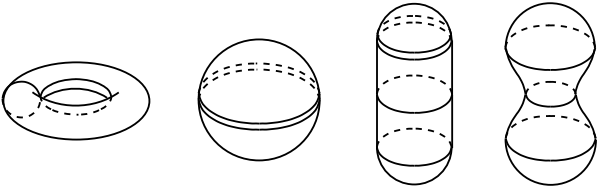}
	\caption{A torus, a sphere, a capsule, and a peanut.}
	\label{fig:manifolds}
\end{figure}

\begin{example}[label=exa:cont]
The torus is multiply connected. A curve that wraps around the torus tube cannot be continuously deformed to a point. A circle around the tube of the torus is a local minimizer of $l$ over the space of closed curves. The sphere $\mathcal{S}^2$ is simply connected. The closed geodesics on $\mathcal{S}^2$ are great circles, \eg the equator. The equator is not a local minimizer of $l$ since there are closed curves of constant latitude arbitrarily close to the equator that are shorter than it, see Fig. \ref{fig:manifolds}. On the capsule, in the regions where the cylinder and hemispheres meet, there are curves which are saddle points of $l$. They are minimizers of $l$ on the cylinder but not on the hemispheres. On both the torus and the peanut there is a single closed curve which is a strict local minimizer of $l$. The torus is multi connected whereas the peanut is simply connected. The condition of simple connectedness is necessary to rule out the existence of a curve of minimum length on a closed manifold, but it is not sufficient.
\end{example}

Assume that the manifold is geodesically complete which implies that there exists at least one geodesic path between any two points $x,y\in\M$. Moreover, assume that for some neigborhood of $x$,  $\mathcal{B}_\varepsilon(x)=\{z\in\M\,|\,d_g(x,z)<\varepsilon\}$, there exists a unique geodesic from $x$ to each $y\in\mathcal{B}_\varepsilon(x)$. The largest value $\varepsilon>0$ for which this holds is the injectivity radius $r(x)$ at $x$. We assume that $\inf_{x\in\M}r(x)>0$. 

The results of this paper concerns the local behaviour of a multi-agent system where the distance $d_g(x_i,x_j)$ between any pair of interacting agents $x_i,x_j\in\M$ can be made arbitrarily small by increasing the number $N$. As such, we can assume we are working on a subset of the manifold where all geodesics are uniquely defined. Under these circumstances, define the notion of a closed broken geodesic:

\begin{definition}
By a closed broken geodesic interpolating a tuple of points $(x_i)_{i=1}^N\in\M^N$ we refer to the closed curve 
\begin{align*}
\gamma=\cup_{i=1}^N \gamma_i\subset\M,
\end{align*}
where $\gamma_i$ is the unique geodesic from $x_i$ to $x_{i+1}$ (using the convention $x_{N+1}=x_1$).
\end{definition}

Given a point $x\in\M$ and a tangent vector $v\in\ts[\M]{x}$, the exponential map $\exp_x:\ts[\M]{x}\rightarrow\M$ yields the point $y\in\M$ that lies at a distance $d_g(v,v)$ from $x$ along the geodesic that passes through $x$ with $v$ as a tangent vector. Let $\mathcal{S}_x\subset\ts[\M]{x}$ be an open set on which $\exp_x$ is a diffeomorphism. Denote $\mathcal{N}_x=\exp_x(\mathcal{S}_x)\subset\M$. The inverse of the exponential map is well-defined on $\mathcal{N}_x$. This inverse is the logarithm map  $\smash{\log_x}:\mathcal{N}_x\rightarrow\ts[\M]{x}$ given by $\log_x:\exp_x(v)\mapsto v$. Define the injectivity radius at $x$, $\inj_x\M$, as the radius of the largest geodesic ball contained in $\mathcal{N}_x$. In the following we assume that $\inf_{x\in\M}\inj_x\mathcal{M}>0$.

The directional derivative of a smooth function $f:\M\rightarrow\R$ at $x\in\M$ along $v\in\ts[\M]{x}$ is given by $\tfrac{\diff}{\diff t}f(\gamma(t))|_{t=0}$, where $\gamma\in\Gamma$ satisfies $\gamma(0)=x$, $\dot{\gamma}(0)=v$. The intrinsic gradient of $f$ is uniquely defined as the vector $\nabla_x f(x)\in\ts[\M]{x}$ which satisfies
\begin{align*}
g_x(\nabla_x f(x),v)=\tfrac{\diff}{\diff t}f(\gamma(t))|_{t=0}
\end{align*}
for all $v\in\ts[\M]{x}$. In particular, $\nabla_x d_g(x,y)^2=-2\log_x(y)$ for all $y\in\mathcal{N}_x$.

A key challenge for control on nonlinear manifolds is to achieve good performance on a global level. It is not possible to achieve global asymptotical stability on a compact manifold by means of continuous, time-invariant feedback \citep{bhat2000topological}. However, it is possible to render an equilibrium or equilibrium set almost globally asymptotically stable (\AGAS{}):
\begin{definition}
An equilibrium set $\mathcal{Q}$ of a dynamical system $\Sigma$ on a Riemannian manifold $(\mathcal{M},g)$ is referred to as almost globally asymptotically stable (\AGAS{}) if it is stable and the flow $\Phi(t,x_0)$ of $\Sigma$ satisfies $\lim_{t\rightarrow\infty}d_g(\mathcal{Q},\Phi(t,x_0))=0$ for all $x_0\in\mathcal{M}\backslash\mathcal{N}$, where the Riemannian measure of the set $\mathcal{N}\subset\mathcal{M}$ is zero.
\end{definition}

\section{Distributed control design on Riemannian manifolds}

We first consider intrinsic consensus, \ie multi-agent consensus on nonlinear spaces described in terms of the intrinsic language of Riemannian geometry. After that we turn to extrinsic consensus, \ie a setting where the Riemannian manifold $(\M,g)$ is embedded in an ambient Euclidean space $\R^{n\times m}$. Averages of agent states' are computed in the ambient space and then projected on the tangent space of the manifold. Some subtle differences in the mathematical challenges posed by each of these two perspectives are pointed out at the end of Section \ref{sec:extrinsic_intro}.

\subsection{Intrinsic consensus}

\noindent Consider a network of $N$ interacting agents. The interaction topology is modeled by a graph $\mathcal{G}=(\mathcal{V},\mathcal{E})$ where the nodes $\V=\{1,\ldots,N\}$ represent agents and an edge $e=\{i,j\}\in\E$ indicates that agent $i$ and $j$ can communicate. Assume that the graph is connected, whereby there is at least an indirect path of communication between any two agents. In this paper, we focus on the cycle graph
\begin{align}\label{eq:HN}
\mathcal{H}_N=(\V,\E)=(\{1,\ldots,N\},\{\{i,i+1\}\,|\,i\in\V\}).
\end{align}
For notational convenience we let $N+1=1$ when adding the indices of $\mathcal{H}_N$, \ie $\{1,N\}\in\E$.

The state $x_i$ of agent $i$ belongs to the manifold $\M$. The states are grouped together in a tuple, $x=(x_i)_{i=1}^N\in\M^N$.

The consensus manifold $\mathcal{C}$ of a  Riemannian manifold $(\mathcal{M},g)$ is the diagonal space of $\M^N$ given by the set
\begin{align}
\mathcal{C}&=\{(x)_{i=1}^N\in\mathcal{M}^N\}.\label{eq:C}
\end{align}
The consensus set is a Riemannian manifold; in fact, it is diffeomorphic to $\mathcal{M}$ by the map $\smash{
(x)_{i=1}^N}\mapsto x$. For any connected graph $\mathcal{G}$, an equivalent definition is
\begin{align*}
\mathcal{C}&=\{(x_i)_{i=1}^N\in\mathcal{M}^N\,|\,x_i=x_j,\,\forall\,\{i,j\}\in\E\}. 
\end{align*}
The terms synchronization and consensus are  interchangeably in this paper. For the Kuramoto model, our notion of consensus is refereed to as phase synchronization \citep{dorfler2014synchronization}.

The states form a dynamical system whose solutions are continuous functions of time for any fixed initial condition. We are interested in the asymptotical behaviour of this system:

\begin{definition}
The agents are said to synchronize, or equivalently, to reach consensus, if $\lim_{t\rightarrow\infty}d_g(x(t),\mathcal{C})=0$. The agents are said to synchronize almost globally, or equivalently, to reach consensus almost globally, if $\mathcal{C}$ is \AGAS.
\end{definition}

Given a graph $(\V,\E)$, define the disagreement function $V:\smash{\mathcal{M}^N}\rightarrow\R$ by
\begin{align}\label{eq:V}
V=\tfrac12\!\!\sum_{\{i,j\}\in\E}\!\!w_{ij}d_g(x_i,x_j)^2,
\end{align}
where $w_{ij}\in(0,\infty)$ and $w_{ji}=w_{ij}$ for all $\{i,j\}\in\E$. The consensus seeking system on $\mathcal{M}$ obtained from $V$ is the gradient descent flow
\begin{align}
\dot{x}=-\nabla V, \quad (\dot{x}_i)_{i=1}^N=(-\nabla_i V)_{i=1}^N,\label{eq:descent}
\end{align}
where $\nabla_i$ denotes the gradient with respect to $x_i$ and $x_i(0)\in\mathcal{M}$. If $\mathcal{G}$ is connected, then by \eqref{eq:C}, $x\in\mathcal{C}$ if and only if $V=0$.

Agent $i$ does not have access to $V$, but can calculate
\begin{align*}
V_i=\tfrac12\sum_{j\in\Ni}w_{ij}d_g(x_i,x_j)^2
\end{align*}
at its current position. Symmetry of $d_g$ gives $V=\tfrac12\sum_{i\in\V}V_i$ whereby it follows that $\nabla_i V_i=\nabla_i V$. From a control design perspective, we can assume that the dynamics of each agent takes the form $\dot{x}_i=u_i$ with $u_i\in\ts[\mathcal{M}]{x_i}$. Since agent $i$ can evaluate $V_i$ at its current position, it is reasonable to assume that it can also calculate $u_i=-\nabla_i V_i$:
\begin{myalgorithm}\label{algo:intrinsic}
	The closed-loop system gradient descent flow of $V$ for $\dot{x}_i=u_i$ under the feedback $u_i=-\nabla_i V_i$ is
	\begin{align}\label{eq:dxidt}
	\dot{x}_i=\sum_{j\in\Ni}w_{ij}\log_{x_i}(x_j),\,\forall\,i\in\V.
	\end{align}
	%*-§
\end{myalgorithm}
The discrete-time equivalent of Algorithm \ref{algo:intrinsic} is introduced by \citet{tron2013riemannian}. Potential shaping is another gradient descent flow approach to consensus; an intrinsic, discrete-time version has been applied to $\SOT$ \citep{tron2012intrinsic}. While potential shaping makes the consensus manifold $\AGAS$, it requires each agent to know some global properties of the graph topology.

\subsection{Extrinsic Consensus}

\label{sec:extrinsic_intro}

\noindent Let the manifold $(\M,g)$ be embedded in an ambient Euclidean space $\R^{n\times m}$. Denote the system state by $\ma{X}=(\ve[i]{X})_{i=1}^N\in\smash{(\R^{n\times m})^N}$. Introduce an extrinsic disagreement function $U:\M\rightarrow[0,\infty)$ given by
\begin{align}
U=\tfrac12\sum_{\{i,j\}\in\E}w_{ij}\|\ve[i]{X}-\ve[j]{X}\|^2.\label{eq:U}
\end{align}
Let $U=\tfrac12\sum_{i\in\V}U_i$, where $U_i=\tfrac12\sum_{j\in \Ni}w_{ij}\|\ve[i]{X}-\ve[j]{X}\|^2$. Note that while the intrinsic disagreement function $V$ is based on geodesic distances, $U$ is based on chordal distances.

Just like for the intrinsic algorithm, $\nabla_i U=\nabla_i U_i$. To calculate the gradient, take any smooth extension $W:\smash{(\R^{n\times m})^N}\rightarrow\R$ of $U$, \ie $W|_{\M}=U$, and utilize that
\begin{align*}
\nabla_i U=\Pi_i\tfrac{\partial}{\partial \ma[i]{X}}W_i,
\end{align*}
where $\nabla_i$ denotes the intrinsic gradient with respect to $\ma[i]{X}$, $\Pi_i:\R^{n\times m}\rightarrow\ts[\M]{\ma[i]{X}}$ is an orthogonal projection map onto $\ts[\M]{\ma[i]{X}}$, $\smash{\tfrac{\partial}{\partial \ma[i]{X}}}$ is the extrinsic gradient in the ambient Euclidean space with respect to $\ma[i]{X}$, and $W=\tfrac12\sum_{i\in\V}W_i$ with $W_i:\R^{n\times m}\rightarrow\R$ being any smooth extension of $U_i$.

A gradient descent flow of $U$ is given by:
\begin{myalgorithm}\label{algo:extrinsic}
	The closed-loop system gradient descent flow of $U$ for $\md[i]{X}=\ma[i]{U}$ under the feedback $\ma[i]{U}=-\nabla_i U_i$ is
	\begin{align}\label{eq:dXidt}
	\md[i]{X}=-\Pi_i\sum_{j\in\Ni}w_{ij}(\ma[i]{X}-\ma[j]{X}),\,\forall\,i\in\V.
	\end{align}
\end{myalgorithm}

\noindent We note that if all $\ma{X}\in\M$ have constant norm, then $\Pi\ma{X}=\ma{0}$, whereby the dynamics \eqref{eq:dXidt} simplifies to
\begin{align}
\md[i]{X}=\Pi_i\sum_{j\in\Ni}w_{ij}\ma[j]{X}.\label{eq:kuramoto}
\end{align}
The assumption of all points having constant Fronbenius norm $r$ implies that the manifold can be embedded in $\R^{n\times m}$ as a subset of a sphere $\mathcal{S}^{nm-1}$ with radius $r$. Examples of such manifolds are the Stiefel and Grassmannian manifolds. Examples of systems on the form \eqref{eq:kuramoto} include, the Kuramoto model, the Lohe model on the $n$-sphere \citep{lohe2010quantum,markdahl2018tac}, the Lohe model on $\SO$ \citep{deville2018synchronization}, and a high-dimensional Kuramoto model on the Stiefel manifold \citep{markdahl2018tac}. The Lohe model on $\U$ \citep{lohe2010quantum,ha2018relaxation} can also be represented by \eqref{eq:kuramoto} via the embedding of $\U$ in $\mathsf{SO}(2n)$.

\begin{remark}
Extrinsic consensus is more important than intrinsic consensus since it is a generalization of the Kuramoto and Lohe models of emergent behaviour in nature. From a technical viewpoint, each case brings its own advantages and challenges. For intrinsic consensus, it is straightforward to relate the curve of minimum length to the disagreement function $V$. However, $V$ is not an analytic function in general, making it more difficult to draw conclusions about the asymptotic behaviour of the system. For extrinsic consensus, it is more difficult to relate the curve of minimum length to the disagreement function $U$. However, $U$ is an analytic function which implies that the gradient descent flow is guaranteed to be well-behaved.
\end{remark}

\subsection{Problem description}

\noindent The aim of this paper is to show that the consensus manifold $\mathcal{C}$ is not an \AGAS{} equilibrium set of (i) the intrinsic gradient descent flow \eqref{eq:dxidt} generated by Algorithm \ref{algo:intrinsic} if $\M$ contains a closed curve $\gamma$ which is a local minimizer of $l$ given by \eqref{eq:l} and (ii) the extrinsic gradient descent flow \eqref{eq:dXidt} generated by Algorithm \ref{algo:extrinsic} if $\M$ is multiply connected. Note that the condition of objective (i) is satisfied by any closed Riemannian manifold that is multiply connected by Theorem \ref{th:simply}.

\subsection{Consensus optimization on Riemannian manifolds}

\noindent Since the system \eqref{eq:dxidt} is a gradient descent flow, it is advantageous to study it from an optimization perspective. %To that end we cover some known results.

\begin{definition}\label{def:minimizer}
A path-connected set $\mathcal{S}\subset\mathcal{X}$ of minimizers of a function $f:\mathcal{X}\rightarrow\R$ is said to be a local minimizer if for some $\delta>0$ there is a ball $\mathcal{B}_\delta(\mathcal{S})=\{x\in\mathcal{X}\,|\,d_g(x,\mathcal{S})<\delta\}$ such that $f|_{\mathcal{S}}\leq f(x)$ for all $x\in\mathcal{B}_\delta(\mathcal{S})$. Moreover, if there is a $\delta>0$ such that the inequality is strict for all $x\in\mathcal{B}_\delta(\mathcal{S})\backslash\mathcal{S}$, then $\mathcal{S}$ is said to be a strict local minimizer.
\end{definition}

\begin{definition}\label{def:critical}
A path-connected set $\mathcal{S}\subset\mathcal{X}$ of minimizers of a function $f:\mathcal{X}\rightarrow\R$ is said to be isolated critical if for some $\delta>0$ there is a ball  $\mathcal{B}_\delta(\mathcal{S})=\{x\in\mathcal{X}\,|\,d_g(x,\mathcal{S})<\delta\}$ such that $\mathcal{S}$ contains all critical points of $f$ in $\mathcal{B}_\delta(\mathcal{S})$.
\end{definition}

\begin{example}[continues=exa:cont]
On the peanut shaped manifold in Fig. \ref{fig:manifolds} there is a set of closed broken geodesics that is a local minimizer and an isolated critical set of $l$. Consider the pill shaped manifold formed by grafting two hemispheres to a cylinder. There is a maximal path-connected set of closed broken geodesics on the cylinder that is an isolated critical set of $l$. It does not consists entirely of local minimizers since there are shorter closed broken geodesics on the hemispheres. There are smaller sets on the cylinder that are local minimizers, however they are not isolated.
\end{example}

\begin{definition}
A graph $\mathcal{G}$ is said to be $(\mathcal{M},g)$-synchronizing if all minimizers of $V$ belong to $\mathcal{C}$.
\end{definition}

\begin{definition}
A graph $\mathcal{G}$ is said to be  $\mathcal{M}$-synchronizing if all minimizers of $U$ belong to $\mathcal{C}$.
\end{definition}

The concept of $\mathcal{S}^1$-synchronizing graphs was introduced to study the Kuramoto model in complex networks \citep{canale2007gluing}. The concept of $\mathcal{M}$-synchronizing graphs is a generalization thereof.

It is not immediate that $\mathcal{G}$ being $(\mathcal{M},g)$-synchronizing implies that $\mathcal{C}$ is an \AGAS{} equilibrium of \eqref{eq:descent}. Since \eqref{eq:descent} is a gradient descent of $V$, $x$ cannot converge to a maximum of $V$. Morover, any saddle point of $V$ is unstable. However, a set of saddle points may still have a region of attraction with positive measure, in which case $\mathcal{C}$ cannot be \AGAS{}. For extrinsic consensus protocols over specific manifolds such as the $n$-sphere \citep{markdahl2018tac} and the Stiefel manifold \citep{markdahl2018prx}, it can be shown that $\mathcal{G}$ being $\mathcal{M}$-synchronizing implies $\mathcal{C}$ is \AGAS{}. For the purpose of this paper we only require the inverse implication for the intrinsic consensus algorithm:

\begin{proposition}\label{prop:not}
Suppose that $\mathcal{G}$ is not $(\M,g)$-synchronizing and that $V$ is $\smash{C^2}$ in an open neighborhood $\mathcal{X}\subset\M$ of a path-connected set of local minimizers $\mathcal{S}$ that is disjoint from  $\mathcal{C}$. If $\mathcal{S}$ consists of strict local minimizers, then it is stable under the closed loop dynamics \eqref{eq:dxidt} of Algorithm \ref{algo:intrinsic}. In particular, $\mathcal{C}$ is not \AGAS{} over  $\mathcal{G}$. If $\mathcal{S}$ is an isolated set of critical points, then it is asymptotically stable.
\end{proposition}

\begin{pf}
The proof is analogous to that of Lyapunov's theorem using $V-V|_{\mathcal{S}}$ as a Lyapunov function, although it applies to sets rather than an equilibrium point. The implication of $\mathcal{C}$ not being \AGAS{} is immediate since stability of $\mathcal{S}$ implies that there is a $\delta>0$ such that the set $\{x\in\mathcal{M}\,|\,d_g(x,\mathcal{S})<\delta\}$, which has positive Riemannian measure, does not belong to the region of attraction of $\mathcal{C}$. For more details, see Appendix \ref{app:proof}.\hfill\phantom{.}\qed
\end{pf}

\section{Main results}
\label{sec:main}

\subsection{Intrinsic consensus}

\label{sec:intrinsic}

Let all the weights in the disagreement function $V$ given by \eqref{eq:V} be equal. Consider the configuration where the agents are distributed equidistantly over a curve of minimum length. It turns out to be a locally optimal solution. Generalizing to the case of unequal weights, we have the following result:

\begin{theorem}\label{prop:main}
Let $(\M,g)$ be a closed, geodesically complete Riemannian manifold. Suppose $\M$ contains a closed curve $\gamma$ of locally minimum length $L=l(\gamma)$. Let $N\geq 3$ and
\begin{align*}
\mathcal{S}=\Bigl\{(x_i)_{i=1}^N\in\M^N\,\big|\,\gamma \Bigl(L\,\tfrac{\sum_{j=1}^{i-1}w_{jj+1}^{-1}}{\sum_{j=1}^Nw_{jj+1}^{-1}}\Bigr)=x_{i},\forall\,i\in\V\Bigr\},
\end{align*}
where $\gamma$ is parametrized by arc length. Any element of $\mathcal{S}$ is a minimizer of the potential function $V$ over $\mathcal{H}_N$. The graph $\mathcal{H}_N$ is not $(\M,g)$-synchronizing. 

Suppose that $V$ is $C^2$ on an open neighborhood of $\mathcal{S}$ and that $w_{12},\ldots,w_{N1}$ satisfy
\begin{align}
d_g(x_i,x_{i+1})=L\frac{w_{ii+1}^{-1}}{\sum_{j=1}^N w_{jj+1}^{-1}}<\inj_{x_i}\M \label{eq:inj}
\end{align}
for all $i\in\V$. If $\mathcal{S}$ is a strictly local minimizer, then it is Lyapunov stable. If $\mathcal{S}$ is an isolated critical set, then it is asymptotically stable. 
\end{theorem}

\begin{pf}
First we show that the elements of $\mathcal{S}$ are locally optimal solutions to
\begin{align}
\label{eq:intrinsic}\begin{split}
\min V&=\tfrac12\sum_{i=1}^N w_{ii+1}d_g(x_i,x_{i+1})^2\\
\textrm{subject to }x_i&\in\M,\,\forall\,i\in\V.
\end{split}
\end{align}
Consider an initial agent configuration in the vicinity of $\gamma$. Let $c$ denote the closed broken geodesic formed by the agents. %Because of \eqref{eq:inj}, there is a neighborhood of $\gamma$ in the space of closed curves where all closed broken geodesics are continuous functions of $x$. 
The constraint 
\begin{align}
l(c)=\sum_{i=1}^N d_g(x_i,x_{i+1})\geq l(\gamma)= L\label{eq:L}
\end{align}
holds given that $c$ is sufficiently close to $\gamma$.

Add the constraint \eqref{eq:L} to the nonlinear program \eqref{eq:intrinsic}, rewriting it as 
\begin{align}\tag{NLP}
\label{eq:constrained}\begin{split}
\min V&=\tfrac12\sum_{i=1}^N w_{ii+1}d_g(x_i,x_{i+1})^2\\
\textrm{subject to }L&\leq\sum_{i=1}^N d_g(x_i,x_{i+1}),\\
x_i&\in\M,\,\forall\,i\in\V.
\end{split}
\end{align}
There is an open neighborhood around $\gamma$ on which the constraint \eqref{eq:L} is redundant whereby any local solution to \eqref{eq:constrained} is a local solution to \eqref{eq:intrinsic} and vice versa. There is hence no loss of generality in restricting our attention to  \eqref{eq:constrained}.

%, \ie the closed broken geodesic $c(y)$ satisfies
%
%\begin{align*}
%l(c(y))=\sum_{i=1}^N d_g(y_i,y_{i+1})\geq L.
%\end{align*} 
%

Consider a point $y\in\M^N$ at which the constraint \eqref{eq:L} holds. Introduce a quadratic program,
\begin{align}
\tag{QP}
\label{eq:QP}\begin{split}
\min f&=\tfrac12\sum_{i=1}^N w_{ii+1}d_{ii+1}^2\\
\textrm{subject to }l(c(y))&=\sum_{i=1}^N d_{ii+1},\\
d_{ii+1}&\geq0,\quad\forall\,i\in\V.
\end{split}
\end{align}
Note that for each solution $y$ to \eqref{eq:constrained}, there is a problem \eqref{eq:QP} which has a corresponding solution given by the vector
\begin{align*}
\ve{d}(y)=[d_{ii+1}(y)]=[d_g(y_i,y_{i+1})]\in[0,\infty)^N.
\end{align*}
The other solutions to \eqref{eq:QP} are not necessarily related to $y$ or to other solutions of \eqref{eq:constrained}. However, if $c(y)$ is a closed geodesic (\ie not just a closed broken geodesic), then each solution $\ve{d}$ to \eqref{eq:QP} generates a solution $x(\ve{d})$ to \eqref{eq:constrained} with the same objective function value, $V(x(\ve{d}))=f(\ve{d})$. In this case, \eqref{eq:QP} can be interpreted as the problem of optimally partitioning the arc length of $l(c)$ into $N$ parts.

The solution $y$ to \eqref{eq:constrained} has the same objective function value as the solution $\ve{d}(y)$ to \eqref{eq:QP}. We show that $\ve{d}(y)$ is suboptimal to \eqref{eq:QP}. However, since $\gamma$ is a closed geodesic, the optimal solution to \eqref{eq:QP} for $c=\gamma$ allows us to obtain a set of solutions to \eqref{eq:constrained}. This is the set $\mathcal{S}$ in Theorem \ref{prop:main}. Let $g(y):\M^N\rightarrow\R$ denote the objective function  value of the optimal solution to \eqref{eq:QP}. We will show that 
\begin{align}
V(y)=f(\ve{d}(y))\geq g(y)\geq g(x)|_{x\subset\gamma}=V|_{\mathcal{S}},\label{eq:train}
\end{align}
which establishes the optimality of $\mathcal{S}$. So far we have only shown the first relation. The second relation follows by the definition of $g$. It remains to establish the last two relations.

The positivity constraint $d_{ii+1}\geq0$ in \eqref{eq:QP} can be relaxed. To see this, note that if $d_{jj+1}<0$ for some $j\in\V$, then $d_{jj+1}$ is counterproductive towards satisfying the constraint \eqref{eq:L} while also incurring a positive cost. A new solution can be constructed where $d_{jj+1}$ is replaced with $0$ while the values of some other variables which assume positive values are decreased so that \eqref{eq:L} still holds. The objective value of the new solution is strictly better than that of the solution from which it was constructed. By relaxing the positive constraints we obtain the equality constrained quadratic program
\begin{align}
\tag{EQP}
\label{eq:QP2}\begin{split}
\min f&=\tfrac12\sum_{i=1}^N w_{ii+1}d_{ii+1}^2\\
\textrm{subject to }l(c(y))&=\sum_{i=1}^N d_{ii+1},\\
d_{ii+1}&\in\R,\quad\forall\, i\in\V.
\end{split}
\end{align}

Equality constrained quadratic programs can be solved using the Lagrange conditions for optimality
\begin{align*}
\begin{bmatrix}
\ma{H} & \mat{A}\\
\ma{A} & \ma{0}
\end{bmatrix}\begin{bmatrix}
\ve{x}\\
\ve{\lambda}
\end{bmatrix}=\begin{bmatrix}
-\ve{c}\\
\ve{b}
\end{bmatrix},
\end{align*}
where $\ma{H}\in\R^{n\times n}$ is the Hessian matrix, $\ma{A}\in\R^{m\times n}$ is the constraint matrix, $\ve{x}\in\R^n$ are the variables, $\ve{\lambda}\in\R^m$ is the vector of Lagrange multiplier, $\ve{c}\in\R^n$ is the coefficients of the linear term in the objective function and $\ve{b}\in\R^m$ is the right-hand side of the constraints \citep{nocedal1999numerical}. For \eqref{eq:QP2} we get
\begin{align*}
\begin{bmatrix}
\ma{W} & \mat{1}\\
\ma{1} & 0
\end{bmatrix}\begin{bmatrix}
\ve{d}\\
\lambda
\end{bmatrix}=\begin{bmatrix}
\ve{0}\\
l
\end{bmatrix},
\end{align*}
where $\ve{d}\in\R^{N}$ is given by $\ve[i]{d}=d_{ii+1}$, $\ve{1}=[1\,\ldots\,1]\smash{\in\R^{N}}$, $\ma{W}$ with $\ma[ii]{W}=w_{ii+1}$ is diagonal, and $l$ is shorthand for the curve length $l(c(y))$.

Denote 
\begin{align*}
\ma{A}=\begin{bmatrix}
\ma{W} & \vet{1}\\
\ve{1} & 0
\end{bmatrix},\quad
\ma{M}=\ma{W}\inv.
\end{align*}
It can be shown that
\begin{align*}
\ma{A}\inv=\frac{1}{\vphantom{1^1}\ve{1}\ma{M}\vet{1}}\begin{bmatrix}
(\ve{1}\ma{M}\vet{1})\ma{M}-\ma{M}\vet{1}\ve{1}\ma{M} & \ma{M}\vet{1}\\\
\ve{1}\ma{M} & -1
\end{bmatrix},
\end{align*}
from which it follows
\begin{align*}
\begin{bmatrix}
\ve{d}\\
\lambda
\end{bmatrix}&=\frac{1}{\vphantom{1^1}\ve{1}\ma{M}\vet{1}}\begin{bmatrix}
(\ve{1}\ma{M}\vet{1})\ma{M}-\ma{M}\vet{1}\ve{1}\ma{M} & \ma{M}\vet{1}\\\
\ve{1}\ma{M} & -1
\end{bmatrix}\begin{bmatrix}
\ve{0}\\
l
\end{bmatrix}\nonumber\\
&=\frac{l}{\vphantom{1^1}\ve{1}\ma{M}\vet{1}}\begin{bmatrix}
\ma{M}\vet{1}\\
-1
\end{bmatrix}.
\end{align*}
The objective value of \eqref{eq:QP2} is
\begin{align}
\sum_{i=1}^N w_{ii+1}d_{ii+1}&=\vet{d}\ma{W}\ve{d}\nonumber\\
&=\Bigl(\frac{l}{\vphantom{1^1}\ve{1}\ma{M}\vet{1}}\Bigr)^2\ve{1}\mat{M}\ma{W}\ma{M}\vet{1}\nonumber\\
&=\frac{l^2(c(y))}{\sum_{i=1}^N w_{ii+1}^{-1}}.\label{eq:VLP}
\end{align}

Recall that $\gamma$ is a local minimizer of $l$ and that $g(y)$ denotes the optimal value to \eqref{eq:QP}. From \eqref{eq:VLP} it follows that 
\begin{align*}
g(y)=\frac{l^2(c(y))}{\sum_{i=1}^N w_{ii+1}^{-1}}\geq\frac{l^2(\gamma)}{\sum_{i=1}^N w_{ii+1}^{-1}}=g(x)|_{x\in\gamma},
\end{align*}
which is the third relation in \eqref{eq:train}. To obtain the last relation in \eqref{eq:train}, $g(x)|_{x\in\gamma}=V|_{\mathcal{S}}$, note that since $\gamma$ is a closed geodesic, any set of points $\smash{x=(x_i)_{i=1}^N\subset\gamma}$ regenerates $\gamma=c(x)$ as their closed broken geodesic. This property allows us to construct a solution $x(\ve{d})$ to \eqref{eq:constrained} from the solution $\ve{d}$ to \eqref{eq:QP2}. The optimal solution to the problem \eqref{eq:QP2} tells us how to position the agent $(x_i)_{i=1}^N$ on $\gamma$ in terms of the arc length distance of $\gamma$ from some arbitrary reference point. The set of such points is $\mathcal{S}$.

%The problems \eqref{eq:constrained} and \eqref{eq:QP2} have the same objective function value at $y$ and $\ve{d}(\ve{y})$. This value is higher than that of the optimal solution to \eqref{eq:QP2} given by \eqref{eq:VLP}. However, if $c$ is a closed geodesic, and not just a closed broken geodesic, then it is possible to construct a solution to \eqref{eq:constrained} from the optimal solution to \eqref{eq:QP2}. In particular, this is possibly for $c=\gamma$, with all such solutions being given by the elements of $\mathcal{S}$. Noting that $\gamma$ minimizes $l(c)$, it is clear that the distribution of agents over $\gamma$ gives a locally optimal solution.

Consider the problem of stability. While the previous discussion only concerned the optimization problem \eqref{eq:intrinsic}, we now need to make sure that the flow \eqref{eq:dxidt} is well-defined. If $w_{12},\ldots,w_{N1}$ satisfy
\begin{align*}
d_g(x_i,x_{i+1})=L\frac{w_{ii+1}^{-1}}{\sum_{i=1}^N w_{ii+1}^{-1}}<\inj_{x_i}\M,
\end{align*}
then the geodesics are unique and the logarithm map is well-defined. Suppose that $V$ is $\smash{C^3}$. The implications that a strict local minimizer $\mathcal{S}$ is stable and an isolated critical set $\mathcal{S}$ is asymptotically stable follows from Proposition \ref{prop:not}.\hfill\phantom{p}\qed\end{pf}

%\begin{figure}[htb!]
%	\centering	\includegraphics[width=0.3\textwidth]{images/torus3}
%	\caption{The agents are perturbed from an optimal distribution on the inner circle of the torus ($\times$) to a closed, broken geodesic which interpolates the perturbed positions ($\wedge$). The resulting configuration is (most likely) suboptimal, given that the agents are restricted to this curve. There is an optimal distribution of agents on this curve ($\vee$). The objective value of that distribution ($\vee$) is higher than that on the inner circle ($\times$). As such, the distribution ($\wedge$) cannot be better than ($\times$).}
%	\label{fig:proof}
%\end{figure}

%Theorem \ref{prop:main} leads to the following result, which states that the non-existence of a curve of minimum length is a necessary condition for $\mathcal{C}$ to be an \AGAS{} equilibrium of \eqref{eq:dxidt} regardless of graph topology.
%
%\begin{corollary}\label{th:main2}
%Suppose the set $\mathcal{S}$ in Theorem \ref{prop:main} consists of strict local minimizers. Then there exists a dense set of initial conditions, a graph $\mathcal{G}$, and a $N\in\N$ such that the closed loop system \eqref{eq:dxidt} generated by Algorithm \ref{algo:intrinsic} does not converge to the consensus manifold $\mathcal{C}$. Whether the consensus manifold $\mathcal{C}$ is \AGAS{} or not depends on $\mathcal{G}$.
%\end{corollary}

\begin{example}[continues=exa:cont]
Theorem \ref{prop:main} establishes asymptotical stability of certain sets on the torus and peanut shaped manifolds in Fig. \ref{fig:manifolds}. Some subsets on the pill shaped manifold are stable, but Theorem \ref{prop:main} does not capture that. The issue is one of having to make sure that certain pathological cases are excluded, see \eg \citet{absil2006stable}.
\end{example}

\subsection{Extrinsic consensus}

\label{sec:extrinsic}

%Here, we repeat the two result of the previous section but for System \label{eq:dXidt}.

%\noindent For the case of extrinsic consensus algorithm, it is no longer straightforward to relate the geometric and topological properties of a general Riemannian manifold to the disagreement function. The chords that interpolate the agents in the ambient space form a polygon rather than a closed curve. Since this polygon can converge to the curve in the limit $N\rightarrow\infty$, it is however reasonable to expect similar behaviour from the systems \eqref{eq:dxidt} and \eqref{eq:dXidt} when $N\gg1$.  

Extrinsic consensus algorithms on Stiefel manifolds $\St$ are known to converge to $\mathcal{C}$ almost globally for any connected graph provided $p\leq\tfrac23n-1$ \citep{markdahl2018prx}. However, almost global convergence does not hold for the Kuramoto model  $\smash{\mathcal{S}^1}\simeq\mathsf{St}(1,2)$ nor $\SOT\simeq\mathsf{St}(2,3)$. In this section we show that these negative findings extend to any closed manifold which is multiply connected.

\begin{theorem}\label{th:extrinsic}
Let $\M$ be a closed, geodesically complete Riemannian manifold which is multiply connected. Then there exists a  set of initial conditions with strictly positive Riemannian measure, a connected graph $\mathcal{G}$, and a $N\in\N$ such that the closed loop system \eqref{eq:dXidt} generated by Algorithm \ref{algo:extrinsic} does not converge to the consensus manifold $\mathcal{C}$.
\end{theorem}

%Given a closed curve $\gamma\subset\M$ which  is not homeomorphic to a point, it is possible to distribute a large number $N$ of agents on $\gamma$ equidistantly with respect to the chordal distance, \ie $\|\ma[i]{X}-\ma[j]{X}\|=k$ for all $\{i,j\}\in\E$. To see this, position the agents along the curve so that all are equidistant $\|\ma[i]{X}-\ma[j]{X}\|=k$ for all $\{i,j\}\in\E\backslash\{1,N\}$. As we increase $k$ from $0$, $\|\ma[1]{X}-\ma[N]{X}\|$ increases from $0$ and ultimately returns to $0$. 

\begin{pf} Since $\M$ is multiply connected, there exists a closed curve $\gamma_0\subset\M$ which is not homeomorphic to a point. For all $N\in\N$, distribute $N$ agents $\smash{(\ma[i]{X})_{i=1}^N}$ approximately equidistantly, in terms of the chordal distance, approximately over $\gamma_0$. Note that there is a set of points with strictly positive Riemannian measure which satisfy this requirement. Let the weights be uniformly bounded away from zero and infinity, \ie there exists an interval  $I\subset(0,\infty)$ such that $w_{ii+1}\in I$ for all $N\in\N$. The  positioning of the agents makes $\|\ma[i+1]{X}-\ma[i]{X}\|$ scale as $\mathcal{O}(\smash{N\inv})$. Then
\begin{align*}
U=\tfrac12\sum_{i=1}^Nw_{ii+1}\|\ma[i+1]{X}-\ma[i]{X}\|^2\in\mathcal{O}(N\inv)
\end{align*}
whereby $U\rightarrow 0$ as $N\rightarrow\infty$. 
	
Consider the level set $\{(\ma[i]{Y})_{i=1}^N\,|\,U\leq \varepsilon\}$. Since $U\rightarrow 0$ as $N\rightarrow\infty$, for any $\varepsilon>0$ there is a $N$ such that $(\ma[i]{X})_{i=1}^N\in\{(\ma[i]{Y})_{i=1}^N\,|\,U\leq \varepsilon\}$. Let $\mathcal{L}$ denote the connected component of the level set which contains $(\ma[i]{X})_{i=1}^N$. Note that for $(\ma[i]{X}(0))_{i=1}^N\in\mathcal{L}$ it holds that $\|\ma[i+1]{X}(t)-\ma[i]{X}(t)\|\leq (2w_{ii+1}^{-1}\varepsilon)^{\frac12}$ for all $t\in[0,\infty)$ since $U$ is decreasing.
	
Let $r(\ma{X})=\inj_{\ma{X}}\!\M$ denote the injectivity radius of $\M$ at $\ma{X}$.  Consider the geodesic ball $\mathcal{D}=\{\ma{Y}\in\M\,|\,d_g(\ma{Y},\ma{X})<r(\ma{X})\}$. There is an open ball $\mathcal{B}$, of strictly positive Lebesgue measure in the ambient space $\R^{n\times m}$, such that $\mathcal{B}\cap\M\subset\mathcal{D}$. Suppose this is false. Then, no matter how small $\mathcal{B}$, there is an $\ma{Z}\in\mathcal{B}\cap\M$ with $d_g(\ma{Z},\ma{X})>r(\ma{X})$. Form a sequence $(\ma[j]{Z})_{j=1}^\infty\subset\M$, such that $\lim_{j\rightarrow\infty}\ma[j]{Z}=\ma{X}$ but $\lim_{j\rightarrow\infty}d_g(\ma[j]{Z},\ma{X})\geq  r(\ma{X})$. That contradicts the continuity of $d_g$. Hence we may chose $N$ large enough that the geodesic distance $d_g(\ma[j]{X}(t),\ma[i]{X}(t))<\inf_{\ma{X}\in\M}r(\ma{X})$, for all $j\in\Ni$ and all $t\in[0,\infty)$. 
	
The geodesic between two agents can be constructed from the  exponential map. The exponential map is a diffeomorphism for $d_g(\ma[j]{X},\ma[i]{X})<\inf_{\ma{X}\in\M}r(\ma{X})$. The closed broken geodesic $\gamma(t)$ that interpolates the agent positions  at each time point can therefore not change discontinuously. The system evolution hence corresponds to a continuous deformation of $\gamma(t)$. If the system  converges to the consensus manifold, then $\lim_{t\rightarrow\infty}d_g(\gamma(t),\mathcal{C})=0$. This contradicts $\gamma_0=\gamma(0)$ not being homeomorphic to a point.\hfill\phantom{p}\qed\end{pf}

%\begin{theorem}
%Suppose that the assumptions of Theorem \ref{th:extrinsic} hold and that the manifold $\M$ is analytic and compact. Then there exists a stable equilibrium set $\mathcal{S}$ which has an empty intersection with the consensus manifold $\mathcal{C}$.
%\end{theorem}

%\begin{pf}
%Consider the level set $\{(\ma[i]{Y})_{i=1}^N\,|\,U\leq \varepsilon\}$. By Theorem \ref{th:extrinsic} there is an $\varepsilon\in(0, \infty)$ there is a dense set $\mathcal{S}$ such that $x\in\mathcal{S}$ implies $\lim_{t\rightarrow\infty}x(t)\notin\mathcal{C}$. Because $\M$ is compact and analytic, $x$ converges to a point $y$. Because $\M$ is compact, the point $y$ is either a saddle point or a minimum of $U$. If $y$ is a minimum point, then we are done by \citep{absil2006stable}. 

%Suppose $y$ is a saddle point. This appears to be a dead end.

%\end{pf}

\begin{remark}
It is possible to relax the assumption of $\G$ being a cycle  to some extent. The graph must allow for the existence of an initial conditions at a distance from consensus such that $U\rightarrow0$ as $N\rightarrow\infty$. For example, we could use a circulant graph where the node degree $d\in2\N$ satisfies $d\ll N$, vertex $i$ is a neighbor of $i-d,\ldots,i+d$, and the same initial condition as for the cycle graph $\mathcal{H}_N$. It is also possible to use a graph that breaks symmetry. For example, we could glue any connected graph $(\mathcal{U},\mathcal{F})$ to node $1$ of $\mathcal{H}_N$. Confine the agents in $\mathcal{U}$ to a ball around agent $1$ that is small enough that $\sum_{\{i,j\}\in\mathcal{F}}w_{ij}d_g(x_i,x_j)^2\in\mathcal{O}(N^{-1})$.
\end{remark}

%\begin{pf}
%The manifold $\mathsf{O}(n)$ is not connected; it is separated  by the function $\det:\mathsf{O}(n)\rightarrow\{-1,1\}$. Achieving consensus on $\mathsf{O}(n)$ from an initial conditions such that one agent has $\det\ma[i]{S}=1$ while another has $\det\ma[j]{S}=-1$ by means of continuous feedback is impossible. The Haar measure on of a set $\mathcal{S}\in\mathsf{O}(n)$ is invariant under left multiplication by any element of $\mathsf{O}(n)$. If $n$ is odd, then multiplication by $-\ma{I}\in\mathcal{O}(n)$ yields a bijective correspondence between $\SO$ and $\mathsf{O}(n)\backslash\SO$. If $n$ is even, then we change the sign of all rows except one. This can be done by left multiplication by an element of $\mathsf{O}(n)$. The Haar measure of $\SO\subset\mathsf{O}(n)$ hence equals that of $\mathsf{O}(n)\backslash\SO$. The Haar measure of neither  set is not zero. It is hence impossible to achieve almost global consensus.\hfill\phantom{x}\qed
%\end{pf}

\begin{figure}[htb!]
\centering	\includegraphics[width=0.35\textwidth]{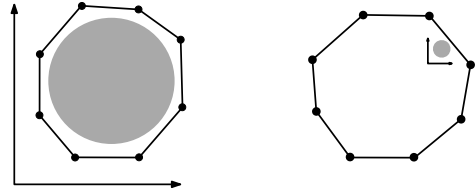}
\caption{\label{fig:extrinsic}The agents ($\cdot$) are approximately equidistantly distributed on a plane with a hole ($\bullet$). Their positions are interpolated by a closed broken geodesic (--).}
\end{figure}

\begin{example}
Consider a manifold consisting of the plane with a hole of radius $r>0$ around a point $\ve{p}\in\R^2$,
\begin{align*}
\M=\{\ve{x}\in\R^2\,|\,\|\ve{x}-\ve{p}\|_2\geq r\}.
\end{align*}
Let the agents be distributed approximately equidistantly in such a manner that the closed broken geodesic $(\ve[i]{x})_{i=1}^N$ encircles the hole. Consider the left case in Fig. \ref{fig:extrinsic}. If the agents are to converge to consensus asymptotically, then there must be a time such that $\|\ve[i+1]{x}-\ve[i]{x}\|_2\geq 2r$ for at least one $i\in\V$. This is not possible for the given initial condition if the weights are equal. By contrast, consider the right case in Fig. \ref{fig:extrinsic}. The agents will pass by the hole without registering that it exists since the chordal distance is not affected. Note that the behaviour of a consensus seeking system is different from that of a curve shortening flow in this respect.

\end{example}

\section{Simulations}

\label{sec:simulation}

Let us use to simulations to explore the question if a manifold being simply connected is not just necessary but also sufficient for the consensus manifold $\mathcal{C}$ of the closed loop system \eqref{eq:dXidt} generated by Algorithm \ref{algo:extrinsic} to be \AGAS{}.

Consider the canonical embedding of the Stiefel manifold $\St$ as a matrix manifold in $\R^{n\times p}$ given by
\begin{align*}
\St=\{\ma{S}\in\R^{n\times p}\,|\,\mat{S}\ma{S}=\ma{I}\}.
\end{align*}
We restate Algorithm \ref{algo:extrinsic} for $\M=\St$ and $w_i=1$:

\begin{myalgorithm}\label{algo:St}
	The input $\ma[i]{U}\in\ts[\St]{i}$ is the negative gradient of the disagreement function, \ie $\ma[i]{U}=-\nabla_i U$. The closed-loop system is a gradient descent flow given by
	\begin{align}
	(\md[i]{S})_{i=1}^N&=-\nabla U=(-\nabla_i U)_{i=1}^N=(-\Pi_i\tfrac{\partial}{\partial \ma[i]{S}} U)_{i=1}^N,\nonumber\\
	\md[i]{S}&=\ma[i]{S}\skews\mat[i]{S}\!\sum_{j\in\Ni}\ma[j]{S}+(\ma{I}\!-\ma[i]{S}\mat[i]{S})\sum_{j\in\Ni}\ma[j]{S},\label{eq:Sid}\\
	\ma[i]{S}(0)&\in\St.
	\end{align}
\end{myalgorithm}
Note that $\skews:\R^{n\times n}\rightarrow \so:\ma{A}\mapsto\tfrac12(\ma{A}-\mat{A})$ is the skew-symmetric part of a matrix. Also note that Algorithm \ref{algo:St} is different from  Algorithm \ref{algo:extrinsic} on $\mathsf{O}(n)$ restricted to $\SO\subset\mathsf{O}(n)$ (see Appendix \ref{app:St}). 

The consensus manifold on $\St$ under the dynamics \eqref{eq:Sid} is \AGAS{} for any connected graph if $p\leq\tfrac23n-1$ \citep{markdahl2018prx}. The results of this paper show that $\mathcal{C}$ cannot be \AGAS{} if the manifold is multiply connected. The manifold $\St$ is simply connected if $p\leq n-2$ \citep{james1976topology}. The only exceptions are $\mathsf{St}(n-1,n)\simeq\SO$, which is multiply connected, and $\mathsf{St}(n,n)\simeq\mathsf{O}(n)$, which is separated by the function $\det:\mathsf{O}(n)\rightarrow\{-1,1\}$. What about the cases of $\tfrac23n-1<p\leq n-2$? We conjecture that $\mathcal{C}$ is \AGAS{} for all such $p$, \ie the that condition $p\leq n-2$ is both necessary and sufficient for $\mathcal{C}$ to be \AGAS{} for all connected graphs.

Let  $\ma{\Phi}:\R\times\St^N\rightarrow\St^N$ denote the flow of \eqref{eq:Sid}, \ie  $\ma{\Phi}(t,(\ma[i,0]{S})_{i=1}^N))=(\ma[i]{S}(t))_{i=1}^N$ given that $(\ma[i]{S}(0))_{i=1}^N=(\ma[i,0]{S})_{i=1}^N\in\St^N$. Denote $\ma{\Phi}=(\ma[i]{\Phi})_{i=1}^N$. Let $\mathcal{R}$ denote the region of attraction of $\mathcal{C}$,
\begin{align*}
\mathcal{R}=\{(\ma[i]{S})_{i=1}^N\in\St^N\,|\,\lim_{t\rightarrow\infty}\ma{\Phi}(t,(\ma[i]{S})_{i=1}^N)\in\mathcal{C}\}.
\end{align*}
The probability measure $\mu(\mathcal{R})$ of $\mathcal{R}$ is the fraction of initial conditions that ultimately yield a consensus. %It is proportional to the size of the `sync basin' \citep{wiley2006size}, \ie the basin of attraction of $\mathcal{C}$.

The probability measure $\mu(\mathcal{R})$ can be calculated by means of Monte Carlo integration:
\begin{align*}
\tfrac1M\sum_{k=1}^M \mathds{1}_{\mathcal{C}}(\lim_{t\rightarrow\infty}\ma{\Phi}(t,(\ma[i]{S}^k)_{i=1}^N))\xrightarrow{\mathrm{a.s.}}\mu(\mathcal{R})\textrm{ as }M\rightarrow\infty,
\end{align*}
where $\mathds{1}:\St^N\rightarrow\{0,1\}$ is the indicator function and $(\ma[i]{S}^k)^N_{i=1}$ for each $k\in\{1,\ldots,M\}$ are samples drawn from the uniform distribution on $\St^N$. To draw a sample $\ma{S}$ from the uniform distribution on $\St$, draw an $\ma{X}\in\R^{n\times p}$ such that each element of $\ma{X}$ is independent and identically normally distributed $N(0,1)$ and form $\ma{S}=\ma{X}(\mat{X}\ma{X})^{-\frac12}$ \citep{chikuse2012statistics}.

A stop criteria is needed to (approximately) calculate $\mathds{1}_{\mathcal{C}}(\lim_{t\rightarrow\infty}\ma{\Phi}(t,(\ma[i]{S})_{i=1}^N))$. Consensus is reached if 
\begin{align}
\max_{\{j,k\}\in\E}\tfrac12\|\ma[j]{\Phi}(T,(\ma[i]{S})_{i=1}^N)-\ma[k]{\Phi}(T,(\ma[i]{S})_{i=1}^N)\|<\varepsilon\label{eq:check}
\end{align}
for a threshold value $\varepsilon\in(0,\infty)$ at a fixed time $T$. %The tensor $\ma{\Phi}(T,(\ma[k,0]{S})_{k=1}^N)$ is obtained by integrating \eqref{eq:Sid} using the function \texttt{ode45} in \textsc{MatLab}. If \eqref{eq:check} is satisfied, we count this as a case of convergence to $\mathcal{C}$. If \eqref{eq:check} is not satisfied the simulation time is extended, with a maximum duration of  $S\gg T$ (the purpose of the check at $t=T$ is simply efficiency). There are potential issues with long simulation times causing numerical errors to accumulate so that either the system leaves the Stiefel manifold or is perturbed from  $\mathcal{R}^\mathsf{c}$ into $\mathcal{R}$. We did not detect any such issues. This approach is biased towards underestimating the value of $\mu(\mathcal{R})$ due to $S<\infty$, but the effect is small for large $S$.

The results are displayed in Table \ref{tab:MC}. Note that they are in agreement with our conjecture that $p\leq n-2$ guarantees convergence to the consensus manifold $\mathcal{C}$. The pairs $(p,n)$ that satisfies $\tfrac23n-1<p\leq p-1$ all have $\mu(\mathcal{R})=1$ (marked in bold). The $(p,n)$ pairs with $1$s in Table \ref{tab:MC} which are not marked in bold satisfy the inequality $p\leq \tfrac23n-1$ whereby $\mathcal{C}$ is \AGAS{} \citep{markdahl2018prx}. Failures to reach consensus occur when $p\in\{n-1,n\}$, \ie for the special orthogonal group and the orthogonal group. For the case of $\mathsf{O}(n)$, the probability that all agents belong to the same connected component is $\smash{2^{-4}}\approx0.06$, which explains the numbers on the diagonal where $p=n$.

%We set $\varepsilon=0.01$, $T=10$, $S=100$.
%\begingroup
%\squeezetable
\begin{table}[thb!]
	\centering
	\normalsize
	\caption{Probability measure, $\mu(\mathcal{R})\in[0,1]$, of the region of attraction $\mathcal{R}$ of the consensus manifold $\mathcal{C}$ on $\St$ for a network given by the graph $\mathcal{H}_5$ defined by \eqref{eq:HN}. The calculation of $\mu(\mathcal{R})$ is done by Monte Carlo integration using $\smash{M=10^4}$ samples of the uniform distribution on $\St$ for each pair $(p,n)$. Rows in the table fix $p$, columns fix $n$. Some cell are left empty since $p\leq n$. %The trivial Stiefel manifold $\mathsf{St}(1,1)=\{-1,1\}$ does not have any dynamics. 
		Bold font indicates that $(p,n)$ satisfies $\tfrac23n-1<p\leq n-2$.\label{tab:MC}} 
	\begin{tabular}{ccccccccccc} %0.06
		%\raisebox{0.12cm}{$\mathcal{G}$\phantom{${}^2$}} &
		%\toprule
		& & & & & $n$ & & & & \\
		& & \multicolumn{1}{c}{2} & \multicolumn{1}{c}{3} & \multicolumn{1}{c}{4} & \multicolumn{1}{c}{5} & \multicolumn{1}{c}{6} & \multicolumn{1}{c}{7} & \multicolumn{1}{c}{8} & \multicolumn{1}{c}{9} \\
		\cline{3-10}
		\rule{0pt}{1em}
		%\midrule
		%\cmidrule{3-11}
		%\addlinespace[1ex]
		%	$\mathcal{S}^1$ & $0$ & $y$ & $0$\\
		%\raisebox{0.12cm}{$\mathcal{G}$\phantom{${}^2$}} &
		%\toprule
		%\midrule
		%\cmidrule{3-11}
		%\addlinespace[1ex]
		%	$\mathcal{S}^1$ & $0$ & $y$ & $0$\\
		& 1 & .95 & 1 & 1 & 1 & 1 & 1 & 1 & 1\\
		& 2 & .05 & .92 & \textbf{1} & 1 & 1 & 1 & 1 & 1\\
		& 3 &  & .06 & .92 & \textbf{1} & 1 & 1 & 1 & 1\\
		& 4 &  &  & .05 & .91 & \textbf{1} & \textbf{1} & 1 & 1\\
		$p$ & 5 &  &  &  & .06 & .89 & \textbf{1} & \textbf{1} & 1\\
		& 6 &  &  &  &  & .05 & .90 & \textbf{1} & \textbf{1}\\
		& 7 &  &  &  &  &  & .06 & .90 & \textbf{1}\\
		& 8 &  &  &  &  &  &  & .06 & .90\\
		& 9 &  &  &  & &  &  &  & .06
		%		711 66
		%		711 66
	\end{tabular}
\end{table}

\section{Conclusions}

\noindent Previous research on almost global synchronization has focused on the graph topology and its influence on convergence \citep{sepulchre2011consensus,dorfler2014synchronization}. In particular, this has been the case for the Kuramoto model on complex networks which is multi-stable for many graphs. By contrast, high-dimensional generalizations of the Kuramoto models to certain  matrix manifolds yield almost global synchronization for any connected graph topology  \citep{markdahl2018tac,markdahl2018prx}. However, the reason for this discrepancy was previously unclear. This paper shows that it can be tied to geometric and topological properties of the manifold. If the manifold is multiply connected, as is \eg the case for the circle $\smash{\mathcal{S}^1}$ and $\SO$, then there exists an obstruction to almost global synchronization over certain graphs. Overcoming this obstruction requires \emph{ad hoc} control design \citep{scardovi2007synchronization,sarlette2009consensus,tron2012intrinsic}. However, in the case of a simply connected manifold, \eg the $n$-sphere for $n\in\N$, such advanced control design techniques are not always needed \citep{markdahl2018tac}. Consider the problem of modeling emergent behaviour and synchronization phenomena that can be observed in nature. This paper suggests that simple and multi connectedness are important properties that modeling needs to account for. For synchronization on multi connected manifolds the Kuramoto model is a good choice. For synchronization on simply connected manifolds the Lohe model on the $n$-sphere is preferable over the Kuramoto model.

%If the manifold is multiply connected, as is \eg the case for the circle $\smash{\mathcal{S}^1}$ and $\SO$, then there exists an obstruction to almost global synchronization over certain graphs. Overcoming this obstruction requires \emph{ad hoc} control design \citep{scardovi2007synchronization,sarlette2009consensus,tron2012intrinsic}. However, in the case of a simply connected manifold, such advanced control design techniques are not always needed. Consider the problem of modeling emergent behaviour and synchronization phenomena that can be observed in nature. This paper suggests that simple and multi connectedness are important properties that modeling needs to account for. For synchronization on multi connected manifolds the Kuramoto model is a good choice. For synchronization on simply connected manifolds the Lohe model on the $n$-sphere is preferable over the Kuramoto model.

\bibliographystyle{named}
\bibliography{autosam}

\appendix

\section{Proof of Proposition \ref{prop:not}}
\label{app:proof}

\begin{pf}
Let $W=V-V|_{\mathcal{S}}$. Note that 
\begin{align*}
\dot{W}&=\tfrac12\!\!\sum_{\{i,j\}\in\E}\!\! w_{ij} g_{x_i}(\nabla_i d_g(x_i,x_j)^2,\dot{x}_i)\\
&=-\!\!\sum_{\{i,j\}\in\E}\!\!g_{x_i}(w_{ij}\log_{x_i}(x_j),\dot{x}_i)\\
&=-\tfrac12\sum_{i\in\V}\sum_{j\in\Ni}g_{x_i}(w_{ij}\log_{x_i}(x_j),\dot{x}_i)\\
&=-\tfrac12\sum_{i\in\V}g_{x_i}(\dot{x}_i,\dot{x}_i),
\end{align*}
which is weakly negative. If (i) $\mathcal{S}$ is a strict local minimizer set (see Definition \ref{def:minimizer}), then there is an open superset $\mathcal{Y}$ of $\mathcal{S}$, $\mathcal{Y}\subset\mathcal{X}$, on which $W$ is strictly positive for all $y\in\mathcal{Y}\backslash\mathcal{S}$. If, moreover, (ii) $\mathcal{S}$ is an isolated critical set (see Definition \ref{def:critical}), then there is a an open superset $\mathcal{Z}$ of $\mathcal{S}$, $\mathcal{Z}\subset\mathcal{Y}$, on which $\dot{W}$ is strictly negative for all $z\in\mathcal{Z}\backslash\mathcal{S}$.

The conditions of Lyapunov's theorem for stability and asymptotic stability are satisfied by (i) and (ii) respectively, except that we are dealing with set stability and intrinsic geometry. The rest of the proof proceeds like the proof of Lyapunov's theorem in \citet{khalil2002nonlinear}, making adjustments for our particular situation as required.

Consider the $(\varepsilon,\delta)$-definition of Lyapunov stability. Given an $\varepsilon>0$, chose an $r\in(0,\varepsilon)$ such that
\begin{align*}\mathcal{B}_r(\mathcal{S})=\{y\in\mathcal{X}\,|\,d_g(y,\mathcal{S})\leq r\}\subset\mathcal{Y}.
\end{align*}
Denote $\alpha=\min_{d_g(y,\mathcal{S})=r}W(y)$. Then $\alpha>0$ because $W$ is strictly positive on $\mathcal{Y}\backslash\mathcal{S}$. Take $\beta\in(0,\alpha)$ and let 
\begin{align*}
\Gamma_\beta=\{x\in \mathcal{B}_r\,|\,W(x)\leq \beta\}.
\end{align*}
Note that by the construction of $\alpha$, $\Gamma_\beta$ cannot contain any points on the boundary of $\mathcal{B}_r(\mathcal{S})$. It follows that $\Gamma_\beta$ is in the interior of $\mathcal{B}_r(\mathcal{S})$. 

For $x(0)\in\Gamma_\beta$, the differential inequality $\dot{W}\leq0$ yields 
\begin{align*}
W(t)\leq W(0)\leq \beta,\quad\forall\,t\in[0,\infty).
\end{align*}
This is interpreted as saying that any trajectory which starts in $\Gamma_\beta$ stays in $\Gamma_\beta$. Existence of a solution $x(t)$ to the gradient descent flow \eqref{eq:descent} for all $t\in\R$ follows from the Picard-Lindelöf theorem \citep{jost2008riemannian}.

Since $W$ is continuous and $0$ on $\mathcal{S}$ it follows that there is an $\delta>0$ such that $d_g(x,\mathcal{S})\leq\delta$ implies that $W<\beta$. Hence there is a closed ball  $\mathcal{B}_\delta(\mathcal{S})\subset\Gamma_\beta$.

The inclusions 
\begin{align*}
\mathcal{B}_\delta(\mathcal{S})\subset\Gamma_\beta\subset\mathcal{B}_r(\mathcal{S})\subset\mathcal{B}_\varepsilon(\mathcal{S})
\end{align*}
imply that $\mathcal{S}$ is stable. This follows since $x(0)\in\mathcal{B}_\delta(\mathcal{S})$ gives $d_g(x(0),\mathcal{S})<\delta$ whereby $x(t)\in\Gamma_\beta$ for all $t\in[0,\infty)$ which in turn yields $d_g(x(t),\mathcal{S})<r<\varepsilon$ for all $t\in[0,\infty)$.

The set $\mathcal{B}_\delta(\mathcal{S})$ has positive Riemannian measure yet it does not belong to the region of attraction of the consensus manifold $\mathcal{C}$. Hence $\mathcal{C}$ is not \AGAS.

The set $\mathcal{S}$ is asymptotically stable if for every $ \nu>0$ there is a $T>0$ such that $d_g(x,\mathcal{S})<\nu$ for all $t\in[T,\infty)$. 

Since $W$ is decreasing and bounded below, $W$ converges as time goes to infinity by the monotone convergence theorem. Denote $c=\lim_{t\rightarrow\infty}d_g(x,\mathcal{S})$. To show that $c=0$, assume that $c>0$ and find a contradiction. Since $W$ is continuous and $W|_{\mathcal{S}}=0$ there is a $\mu>0$ such that 
\begin{align*}
\mathcal{B}_\mu(\mathcal{S})\subset\{x\in\mathcal{B}_r\,|\,W< c\}\cap\mathcal{Z}.
\end{align*}
The trajectory $x(t)$ cannot enter $\mathcal{B}_\mu(\mathcal{S})$ for any $t\in[0,\infty)$ because that would result in $W(x(s))<c$ for all $s\geq t$. The trajectory is hence confined to the set 
\begin{align*}
\mathcal{A}=\{x\in\mathcal{Z}\,|\,\mu\leq d_g(x,\mathcal{S})\leq r\}.
\end{align*}
The minimum rate of decrease on $\mathcal{A}$, $\zeta=\min_{x\in\mathcal{A}}|\dot{V}|$, exists by virtue of the Weierstrass extreme value theorem since $W$ is continuous and the set $\mathcal{A}$ is bounded. Note that $\zeta>0$ since $\mathcal{S}$ is an isolated set of critical points. From
\begin{align*}
W(x(t))=W(x(0))+\int_0^t \dot{W}(x(\tau))\diff\tau\leq W(x(0))-\zeta t
\end{align*}
it follows that there is a time $s$ such that $W(x(s))=0$. This contradicts the assumption of $\lim_{t\rightarrow\infty}W(t)=c>0$.\hfill\phantom{p}\qed\end{pf}

\section{A consensus protocol on $\mathsf{St}(n-1,n)$}
\label{app:St}

Consider the canonical embedding of $\mathsf{O}(n)$ as a matrix manifold in the ambient space $\smash{\R^{n\times n}}$. Let $\smash{(\ma[i]{Q})_{i=1}^N}$ denote the state of a multi-agent system on $\mathsf{O}(n)$. The disagreement function $U$ can be simplified as
\begin{align*}
U&=\sum_{e\in\E}\|\ve[i]{Q}-\ma[j]{Q}\|^2=2\sum_{e\in\E}n-\langle\ma[j]{Q},\ve[i]{Q}\rangle,
\end{align*}
where a factor of $2$ has been introduced for notational convenience. Algorithm \ref{algo:extrinsic} in the special case of $\M=\mathsf{O}(n)$ is:
\begin{myalgorithm}\label{algo:SO}
The input $\ma[i]{U}\in\ts[\mathsf{O}(n)]{i}$ is the negative gradient of the disagreement function, \ie $\ma[i]{U}=-\nabla_i U$. The closed-loop system is a gradient descent flow given by
\begin{align}
(\md[i]{Q})_{i=1}^N&=-\nabla U=(-\nabla_i U)_{i=1}^N=(-\Pi_i\tfrac{\partial}{\partial \ma[i]{Q}} U)_{i=1}^N\nonumber\\
\md[i]{Q}&=2\ma[i]{Q}\skews\mat[i]{Q}\sum_{j\in\Ni}\ma[j]{Q}\nonumber\\
&=\sum_{j\in\Ni}\ma[j]{Q}-\ma[i]{Q}\mat[j]{Q}\ma[i]{Q}.\label{eq:Rid}
\end{align}
\end{myalgorithm}

For the restriction to $\SO\subset\mathsf{O}(n)$ and a circulant graph $\mathcal{G}$, there are stable equilibria $\smash{(\ma[i]{R})_{i=1}^N\notin\mathcal{C}}$ \citep{deville2018synchronization}.

Consider the relation between Algorithm \ref{algo:SO} on $\SO$ and Algorithm \ref{algo:St} on $\St$ in the case of $p=n-1$ for which $\St\simeq\SO$. 
	
Let $\ma[i]{T}$, $\ve[i]{t}$ be such that $\ma[i]{R}=[\ma[i]{T}\,\ve[i]{t}]\in\SO$ where $\ma[i]{R}$ follow the dynamics generated by Algorithm \ref{algo:SO} on $\SO$. Then
\begin{align*}
\md[i]{R}&=[\md[i]{T}\,\vd[i]{t}]=\sum_{j\in\Ni}[\ma[j]{T}\,\ve[j]{t}]-[\ma[i]{T}\,\ve[i]{t}]\begin{bmatrix}\mat[j]{T}\\
\vphantom{1^{1^1}}\vet[j]{t}\end{bmatrix}[\ma[i]{T}\,\ve[i]{t}]\\
&=\sum_{j\in\Ni}[\ma[j]{T}\,\ve[j]{t}]-(\ma[i]{T}\mat[j]{T}+\ve[i]{t}\vet[j]{t})[\ma[i]{T}\, \ve[i]{t}]
\end{align*}
which yields
\begin{align}
\md[i]{T}&=\sum_{j\in\Ni}\ma[j]{T}-\ma[i]{T}\mat[j]{T}\ma[i]{T}-\ve[i]{t}\vet[j]{t}\ma[i]{T}\nonumber\\
&=2\ma[i]{T}\skews\mat[i]{T}\!\sum_{j\in\Ni}\!\ma[j]{T}+2\ve[i]{t}\vet[i]{t}\!\sum_{j\in\Ni}(\ma[j]{T}-\tfrac12\ve[i]{t}\vet[j]{t}\ma[i]{T}),\nonumber\\
\vd[i]{t}&=\sum_{j\in\Ni}\ma[j]{t}-\ma[i]{T}\mat[j]{T}\ve[i]{t}-\ve[i]{t}\vet[j]{t}\ve[i]{t}\nonumber\\
&=2\skews\bigl(\sum_{j\in\Ni}\ma[j]{T}\mat[i]{T}\bigr)\ve[i]{t}+\sum_{j\in\Ni}\ma[j]{t}-\langle\ve[j]{t},\ve[i]{t}\rangle\ve[i]{t}.\label{eq:SO}
\end{align}

Let $\ma[i]{S}\in\mathsf{St}(n-1,n)$ follow the dynamics of Algorithm \ref{algo:St}. Let $\ve[i]{s}\in\R^n$ be such that $[\ma[i]{S}\,\ve[i]{s}]\in\SO$. Then $\mat[i]{S}\ve[i]{s}=\ve{0}$, \ie $\md[i]{S}{}\!\!\mtr\ve[i]{s}+\mat[i]{S}\vd[i]{s}=\ve{0}$. This yields $\mat[i]{S}\vd[i]{s}=-\sum_{j\in\Ni}\mat[j]{S}\ve[i]{s}$. Since $\ma[i]{S}\mat[i]{S}+\ve[i]{s}\vet[i]{s}=\ma{I}$ and $\vet[i]{s}\vd[i]{s}=0$, we obtain $\vd[i]{s}=-\ma[i]{S}\sum_{j\in\Ni}\mat[j]{S}\ve[i]{s}=\skews(\sum_{j\in\Ni}\ma[j]{S}\mat[i]{S})\ve[i]{s}$. The system \eqref{eq:Sid} on $\mathsf{St}(n-1,n)$ hence yields the following system on $\SO$:
\begin{align}
\label{eq:St}\begin{split}
\md[i]{S}&=\ma[i]{S}\skews\mat[i]{S}\sum_{j\in\Ni}\ma[j]{S}+\ve[i]{s}\vet[i]{s}\sum_{j\in\Ni}\ma[j]{S},\\
\vd[i]{s}&=\skews\bigl(\sum_{j\in\Ni}\ma[j]{S}\mat[i]{S}\bigr)\ve[i]{s}.
\end{split}
\end{align}

The systems \eqref{eq:SO} and \eqref{eq:St} yield different paths on $\SO$ as is clear by inspection.

\end{document}

%% file: myIncludes.tex
\usepackage{subdepth} % fixes issues with superscript and subscript vertical position when both are used (there is still some issue with horizontal placement, look at $\|x\|_2^2$)

\usepackage{graphics} % for pdf, bitmapped graphics files
\usepackage{epsfig} % for postscript graphics files

\usepackage{amssymb} % amsfonts % assumes amsmath package installed
\usepackage{amsmath,psfrag,color,graphicx} %euscript

\usepackage{thmtools}
\renewcommand\thmcontinues[1]{Continued}

\usepackage{graphicx}

\usepackage[small,hang,bf,up,labelsep=period]{caption} % have not checked that it works, no figures yet
\usepackage{quoting} % for abstracts
\usepackage{lettrine} % for the first letter
\usepackage{tocloft} % for the table of contents title
\usepackage{titletoc}
\usepackage{fmtcount}
\usepackage{calc} % to use arithmics with counters
\usepackage{cite} 
\usepackage{multicol} % this could be used to create two column chapters
\usepackage{centernot}
\usepackage{nicefrac}
\usepackage{dsfont}
\usepackage{booktabs}
\usepackage{mathtools, cuted}

\usepackage{lipsum}
\usepackage{psfrag}
\usepackage[export]{adjustbox} %\adjincludegraphics is more powerful than includegraphics
\PassOptionsToPackage{cmyk}{xcolor}
\usepackage[cmyk]{xcolor}
\selectcolormodel{cmyk}

\usepackage{textcomp}

\usepackage{txfonts}
\usepackage{lmodern}
\usepackage{times}

\usepackage{scalerel}[2016/12/29]
\DeclareFontEncoding{LS1}{}{}
\DeclareFontSubstitution{LS1}{stix}{m}{n}

\makeatletter
\newcommand*\textmathversion{\csname textmv@\math@version\endcsname}
\newcommand*\textmv@normal{m}
\newcommand*\textmv@bold{b}
\makeatother

\newcommand{\ve}[2][]{\ensuremath{\boldsymbol{\mathrm{#2}}}_{#1}}
\newcommand{\vet}[2][]{\ensuremath{\smash{\boldsymbol{\mathrm{#2}}^{\!\top}_{#1}}}}
\newcommand{\vd}[2][]{\ensuremath{\dot{\boldsymbol{\mathrm{#2}}}_{#1}}}

\newcommand{\ma}[2][]{\ensuremath{\boldsymbol{\mathrm{#2}}}_{#1}}
\newcommand{\mat}[2][]{\ensuremath{\boldsymbol{\mathrm{#2}}^{\!\top}_{#1}}}
\newcommand{\md}[2][]{\ensuremath{\dot{\boldsymbol{\mathrm{#2}}}}_{#1}}

\DeclareMathOperator{\skews}{\ensuremath{\mathrm{skew}}}

\newcommand{\R}{\ensuremath{\mathds{R}}}

\newcommand{\N}{\ensuremath{\mathds{N}}}

\newcommand{\SOT}{\ensuremath{\mathsf{SO}(3)}}

\newcommand{\SO}{\ensuremath{\mathsf{SO}(n)}}

\newcommand{\so}{\ensuremath{\mathsf{so}(n)}}

\newcommand{\St}{\ensuremath{\mathsf{St}(p,n)}}

\newcommand{\U}{\ensuremath{\mathsf{U}(n)}}

\DeclareMathOperator{\inj}{\mathrm{inj}}

\newcommand{\ie}{\textit{i.e.}, }

\newcommand{\eg}{\textit{e.g.}, }

\newcommand{\mtr}{\hspace{-0.3mm}\ensuremath{^\top}}

\newcommand{\diff}{\ensuremath{\mathrm{d}}}

\newcommand{\inv}{\ensuremath{^{-1}}}

\newcommand\raiseT[2]{\setbox0\hbox{$#1{#2}$}\raise\dp0\box0}

\newcommand{\G}{\ensuremath{\mathcal{G}}}

\newcommand{\V}{\ensuremath{\mathcal{V}}}

\newcommand{\E}{\ensuremath{\mathcal{E}}}

\newcommand{\M}{\ensuremath{\mathcal{M}}}

\DeclareMathOperator{\trace}{\ensuremath{\mathrm{tr}}}

\newcommand{\Ni}{\ensuremath{\mathcal{N}_i}}

\makeatletter
\newcommand{\raisemath}[1]{\mathpalette{\raisem@th{#1}}}
\newcommand{\raisem@th}[3]{\raisebox{#1}{$#2#3$}}
\makeatother

\newcommand{\ts}[2][]{\ensuremath{\mathsf{T}_{#2}#1}}

\newcommand{\AGAS}{\textsc{agas}}

\definecolor{kthbluergb}{RGB}{25,84,166}%{0,108,183}
\definecolor{kthbluecmyk}{cmyk}{1,0.55,0,0}

\definecolor{kthblueA}{RGB}{25,84,166}%{0,108,183}
\definecolor{kthblueB}{RGB}{46,124,192}%{0,108,183}
\definecolor{kthblueC}{RGB}{112,153,209}%{0,108,183}
\definecolor{kthblueD}{RGB}{164,186,225}%{0,108,183}
\definecolor{kthblueE}{RGB}{211,220,241}%{0,108,183}

\newcounter{counter} % there is only one counter

\newtheorem{theorem}[counter]{Theorem}

\newtheorem{proposition}[counter]{Proposition}

\newtheorem{definition}[counter]{Definition}
\newtheorem{example}[counter]{Example}

\newtheorem{myalgorithm}[counter]{Algorithm}

\newtheorem{remark}[counter]{Remark}


\newcounter{parentnumber}